\documentclass{amsart}

\usepackage{euscript}

\usepackage[dvips]{graphicx}

\newtheorem{thm}{Theorem}
\newtheorem{pr}{Proposition}
\newtheorem{lem}{Lemma}

\newtheorem{prop}{Property}

\theoremstyle{definition}
\newtheorem{defn}{Definition}
\newtheorem{rem}{Remark}
\newtheorem{example}{Example}

\DeclareMathOperator{\Imp}{\mathrm{Im}\,}
 \DeclareMathOperator{\Rep}{\mathrm{Rep}\,}
\DeclareMathOperator{\End}{\mathrm{End}}
\DeclareMathOperator{\Hom}{\mathrm{Hom}}
\DeclareMathOperator{\Idem}{\mathrm{Idem}}

\author{Yu. P. Moskaleva}
\address{Taurida National University, 4 Vernads'ky, Simferopol, 95007, Ukraine}
\email{YulMosk@mail.ru}

\author{Yu. S. Samo\v\i{}lenko}
\address{Institute of Mathematics, National Academy of
  Sciences of Ukraine, 3 Teresh\-chenkivs'ka, Kyiv, 01601,
  Ukraine}
\email{yurii\_sam@imath.kiev.ua}
\thanks{This research was partially supported by the State Foundation for Fundamental Research of Ukraine,
grant no.~01.07/071 and by the DFG, grant no.~436UKR 113/71/0-1.}

\date{22/11/2005}

\title[Systems of $n$ subspaces and representations of $*$-algebras\dots]%
{Systems of $n$ subspaces and representations of  $*$-algebras
  generated by projections}

\begin{document}

\begin{abstract}
  In this paper we study a relationship between systems of $n$
  subspaces and representations of $*$-algebras generated by
  projections. We prove that irreducible nonequivalent
  $*$-representations of $*$-algebras $\mathcal P_{4,com}$ generate
  all nonisomorphic transitive quadruples of subspaces of a finite
  dimensional space.
\end{abstract}

\maketitle

\section{Introduction}

There are many articles that deal with a description of systems
$S=(H;H_1,H_2,\ldots,H_n)$ of $n$ subspaces $H_i$, $i=\overline{1,n}$,
of a Hilbert space $H$, which can be infinite or finite dimensional,
up to an isomorphism or the unitary equivalence.

In particular, transitive quadruples of subspaces (see Section~2) of a
finite dimensional space were described in~\cite{B}, indecomposable
quadruples were found in~\cite{GP,N}.

For a pair of subspaces $H_1$, $H_2$ of a Hilbert space $H$ there is a
structure theorem (see, for example,~\cite{Ha1}) that describes pairs
of orthogonal projections onto these subspaces, up to the unitary
equivalence, in terms of sums or integrals of irreducible one- or
two-dimensional pairs of orthogonal projections. For three subspaces,
to get such a theorem is unrealistic, --- the problem of getting a
unitary description of $n$ orthogonal projections for $n\geq3$ is
$*$-wild (see~\cite{KS1,KS2}). So, if we need to get a description of
collections of $n$ orthogonal projections up to the unitary
equivalence, it is necessary to introduce additional relations. Recent
works of Ukrainian mathematicians (see~\cite{ZS,OS} and the
bibliography therein) make a study of irreducible systems of
orthogonal projections $P_1, P_2, \ldots, P_n$ such that their sum is
a multiple of the identity operator.

In~\cite{EW}, the authors suspect that there is a relationship between
systems of $n$ subspaces and representations of $*$-algebras generated
by projections, --- ``There seems to be interesting relations with the
study of $*$-algebras generated by idempotents by S.~Kruglyak and
Yu.~Samoilenko~\cite{KS2} and the study on sums of projections by
S.~Kruglyak, V.~Rabanovich and Yu.~Samoilenko~\cite{KRS}. But we do
not know the exact implication \dots''~\cite{EW}. This paper is
devoted to a study of this relationship.

For an irreducible collection of orthogonal projections, $P_1$, $P_2$,
$\ldots$, $P_n$, on a Hilbert space $H$ such that $\sum_{i=1}^{n}
P_i=\alpha I_H$, consider the system of $n$ subspaces
$S=(H;P_1H,P_2H,\ldots, P_nH)$. Let us formulate the following
hypothesis: collections of orthogonal projections such that their sum
is a multiple of the identity operator, that is, irreducible
nonequivalent $*$-representations of the $*$-algebras $\EuScript
P_{n,com}$ (see Section~3), generate nonisomorphic transitive systems.
In Section~4, we prove this hypothesis for $n=1$ and $n=2$. There,
irreducible nonequivalent $*$-representations of the $*$-algebras
$\EuScript P_{1,com}$ and $\EuScript P_{2,com}$ generate all
nonisomorphic transitive systems of one or two subspaces in an
arbitrary Hilbert space. We also prove there that, for $n=3$,
irreducible nonequivalent $*$-representations of the $*$-algebra
$\EuScript P_{3,com}$ generate all nonisomorphic transitive systems of
three subspaces of a finite dimensional linear space. Let us remark
that it is an unsolved problem to describe irreducible triples of
subspaces of an infinite dimensional space or even to prove their
existence for $n=3$ (see~\cite{Ha2}). If $n=4$, we prove in Section~4
that ireducible nonequivalent $*$-representations of the $*$-algebras
$\EuScript P_{n,com}$ generate all nonisomorphic transitive systems for
a finite dimensional space. Since irreducible nonequivalent
$*$-representations of the $*$-algebra $\EuScript P_{4,com}$ can only
be finite dimensional, irreducible nonequivalent $*$-representations
of the $*$-algebra $\EuScript P_{4,com}$ already do not generate all
nonisomorphic transitive systems of four subspaces if $n=4$, see, for
example,~\cite{EW} and the bibliography therein.

\section{Systems of $n$ subspaces}
\subsection{Definitions and main properties}

All statements of this section are regarded as known (see, for
example,~\cite{EW,OS}) and given without proofs. Let $H$ be a Hilbert
space, $H_1$, $H_2$, $\ldots$, $H_n$ be $n$ subspaces of the space
$H$. Denote by $S=(H;H_1,H_2,\ldots,H_n)$ the system of $n$ subspaces
of the space $H$. Let $S=(H;H_1,H_2,\ldots,H_n)$ be a system of $n$
subspaces of the Hilbert space $H$ and
$\tilde{S}=(\tilde{H};\tilde{H}_1,\tilde{H}_2,\ldots,\tilde{H}_n)$ a
system of $n$ subspaces of the Hilbert space $\tilde{H}$.

\begin{defn}\label{d1}
  A linear mapping $R:H\rightarrow\tilde{H}$ of the space $H$ into the
  space $\tilde{H}$ is called a homomorphism of the system $S$ into
  the system $\tilde{S}$ and denoted by $R:S\rightarrow\tilde{S}$, if
  $$
  R(H_i)\subset\tilde{H}_i,\quad i=\overline{1,n}.
  $$
\end{defn}

\begin{defn}\label{d2}
  A homomorphism $R:S\rightarrow\tilde{S}$ of a system $S$ into a
  system $\tilde{S}$ is called an isomorphism, and denoted by
  $R:S\rightarrow\tilde{S}$, if the mapping $R:H\rightarrow\tilde{H}$
  is a bijection and $R(H_i)=\tilde{H}_i,\forall
  i=\overline{1,n}$.
\end{defn}

Systems $S$ and $\tilde{S}$ will be called isomorphic and denoted by
$S\cong\tilde{S}$, if there exists an isomorphism
$R:S\rightarrow\tilde{S}$.

\begin{defn}\label{d3}
  We say that systems $S$ and $\tilde{S}$ are unitary equivalent, or
  simply equivalent, if $S\cong\tilde{S}$ and the isomorphism
  $R:S\rightarrow\tilde{S}$ can be chosen as to be a unitary operator.
\end{defn}

For each system $S=(H;H_1,H_2,\ldots,H_n)$ of $n$ subspaces of a
Hilbert space $H$ there is a naturally connected system of orthogonal
projections $P_1$, $P_2$, $\ldots$, $P_n$, where $P_i$ is the
orthogonal projection operator onto the subspace $H_i$,
$i=\overline{1,n}$. A system of projections $P_1$, $P_2$, $\ldots$,
$P_n$ on a Hilbert space $H$ such that $\Imp P_i=H_i$ for
$i=\overline{1,n}$ will be called a system of orthogonal projections
corresponding to the system of subspaces
$S=(H;H_1,H_2,\ldots,H_n)$. And conversely, for each system of
projections there is a naturally connected system of subspaces. The
system $S=(H;P_1H, P_2H,\ldots,
P_nH)$ will be called a system generated by the system of the
projections $P_1$, $P_2$, $\ldots$, $P_n$.

\begin{defn}\label{d4}
  A system of orthogonal projections $P_1$, $P_2$, $\ldots$, $P_n$ on
  a Hilbert space $H$ is called unitary equivalent to a system
  $\tilde{P}_1$, $\tilde{P}_2$, $\ldots$, $\tilde{P}_n$ on a Hilbert
  space $\tilde{H}$ if there exists a unitary operator
  $R:H\rightarrow\tilde{H}$ such that $RP_i=\tilde{P}_iR$,
  $i=\overline{1,n}$.
\end{defn}

It is clear that systems $S$ and $\tilde{S}$ are unitary equivalent if
and only if the corresponding systems of orthogonal projections are
unitary equivalent.

\begin{prop}\label{p1}
  Let $S=(H;H_1,H_2,\ldots,H_n)$,
  $\tilde{S}=(\tilde{H};\tilde{H}_1,\tilde{H}_2,\ldots,\tilde{H}_n)$
  be systems of $n$ subspaces of  Hilbert spaces $H$ and
  $\tilde{H}$. Let $P_i$ and $\tilde{P}_i$ be orthogonal projection
  operators onto $H_i$ and $\tilde{H}_i$, correspondingly, $
  i=\overline{1,n}$. The systems $S$ and $\tilde{S}$ are isomorphic if
  and only if there exists an invertible operator
  $T:H\rightarrow\tilde{H}$ such that
  $$
  P_i=T^{-1}\tilde{P}_iTP_i\quad\mbox{и}\quad
  \tilde{P}_i=TP_iT^{-1}\tilde{P}_i,\quad i=\overline{1,n}.
  $$
\end{prop}

\begin{rem}\label{r1}
  If systems $S$ и $\tilde{S}$ are unitary equivalent, then
  $S\cong\tilde{S}$. The converse is not true.
\end{rem}

Denote by $\Hom(S,\tilde{S})$ the set of homomorphisms of the system
$S$ into the system $\tilde{S}$, and by $\End(S):=\Hom(S,S)$ the algebra
of endomorphisms from $S$ into $S$, that is,
$$
\End(S)=\{R\in B(H)|R(H_i)\subset H_i, i=\overline{1,n}\}.
$$

\begin{defn}\label{d5} 
  A system $S=(H;H_1,H_2,\ldots,H_n)$ of $n$ subspaces of a space $H$
  will be called transitive if $\End(S)=\mathbb C I_H$.
\end{defn}

\begin{rem}
  Isomorphic systems are simultaneously either transitive or
  nontransitive.  
\end{rem}

Let us introduce the notion of an indecomposable system, which is
equivalent to the definition used in~\cite{GP,EW}. Denote
$$
\Idem(S)=\{R\in B(H)|R(H_i)\subset H_i, i=\overline{1,n},R^2=R\}.
$$

\begin{defn}\label{d6}
  A system $S=(H;H_1,H_2,\ldots,H_n)$ of $n$ subspaces  of a space $H$
  will be called indecomposable if $\Idem(S)=\{0,I_H\}$.
\end{defn}

\begin{rem}
  Isomorphic systems are simultaneously decomposable or
  indecomposable.
\end{rem}

\begin{defn}\label{d7}
  A system of orthogonal projections  $P_1$, $P_2$, $\ldots$, $P_n$ on
  a Hilbert space $H$, which possesses only trivial invariant
  subspaces, is called irreducible.
\end{defn}

\begin{rem}
  Systems of unitary equivalent systems of orthogonal projections are
  simultaneously reducible or irreducible.
\end{rem}

The following proposition answers the question about a relation
between the notions of a transitive system, an indecomposable system,
irreducibility of the corresponding system of orthogonal projections.

\begin{pr}\label{u1}
  If a system of subspaces is transitive, then it is indecomposable.
  If a system of subspaces is indecomposable, then the corresponding
  system of orthogonal projections is irreducible.
\end{pr}

\begin{proof}
  The first statement follows from the obvious inclusion
  $\Idem(S)\subset \End(S)$ and the definitions of a transitive and
  an indecomposable systems. To prove the second statement, we use the
  Schur's lemma (see, for example,~\cite{OS}). A system of orthogonal
  projections  $P_1$,  $P_2$, $\ldots$, $P_n$ on a Hilbert space $H$
  is irreducible if and only if $\{R\in
  B(H)|RP_i=P_iR, i=\overline{1,n}, R^2=R, R^*=R\}=\{0,I_H\}$. The
  identity $\{R\in B(H)|RP_i=P_iR, i=\overline{1,n}, R^2=R,
  R^*=R\}=\{R\in B(H)|R(\Imp P_i)\subset\Imp P_i, i=\overline{1,n},
  R^2=R, R^*=R\}$, on the one hand, and the inclusion $\{R\in
  B(H)|R(H_i)\subset H_i, i=\overline{1,n}, R^2=R, R^*=R\}\subset
  \Idem(S)$, on the other hand, finish the proof.
\end{proof}

\begin{example}
  Let $S=(\mathbb C^2;\mathbb C(1,0),\mathbb
  C(\cos\theta,\sin\theta))$, $\theta\in(0,\pi/2)$ and
  $\tilde{S}=(\mathbb C^2;\mathbb C(1,0),\mathbb C(0,1))$. The
  decomposable system $S$, which corresponds to the irreducible pair of
  orthogonal projections, is isomorphic but not unitary equivalent to
  the decomposable system $\tilde{S}$ that corresponds to the reducible
  pair of orthogonal projections.
\end{example}

\begin{defn}\label{d8}
  Let $S=(H;H_1,H_2,\ldots,H_n)$ be a system of $n$ subspaces of a
  Hilbert space $H$. By an orthogonal complement to the system $S$, we
  will call the system
  $S^\perp=(H;H^\perp_1,H^\perp_2,\ldots,H^\perp_n)$.
\end{defn}

\begin{prop}\label{p2}
  Let $S=(H;H_1,H_2,\ldots,H_n)$ be a system of $n$ subspaces of a
    Hilbert space $H$. Then $S$ is transitive (indecomposable) if and
    only if  $S^\perp$ is transitive (indecomposable).
\end{prop}

Property~\ref{p2} follows directly, since if $R:S\rightarrow\tilde{S}$
is a homomorphism of the system $S$ into $\tilde{S}$, then
$R^*:\tilde{S}^\bot\rightarrow S^\bot$ is a homomorphism of the system
$\tilde{S}$ into $S$, because, if $R:H\rightarrow\tilde{H}$ is a
linear operator such that $R(H_i)\subset\tilde{H}_i,\forall
i=\overline{1,n}$, then $R^*:\tilde{H}\rightarrow H$ and
$R^*(\tilde{H}^\bot_i)\subset H^\bot_i,\forall i=\overline{1,n}$.
\begin{defn}\label{d9}
  Let $S=(H;H_1,H_2,\ldots,H_n)$ and
  $\tilde{S}=(\tilde{H};\tilde{H}_1,\tilde{H}_2,\ldots,\tilde{H}_n)$
  be two systems of $n$ subspaces. We say that $S\cong\tilde{S}$ up to
  a rearrangement of subspaces if there is a permutation $\sigma\in
  S_n$ such that the systems $\sigma(S)$ and $\tilde{S}$ are
  isomorphic, where
  $\sigma(S)=(H;H_{\sigma(1)},H_{\sigma(2)},\ldots,H_{\sigma(n)})$,
  that is, there exists and invertible operator
  $R:H\rightarrow\tilde{H}$ such that $R(H_{\sigma(i)})=\tilde{H}_i$,
  $\forall i=\overline{1,n}$.
\end{defn}

\subsection{Transitive systems of one, two, and three subspaces}

In this section we give a description of transitive systems of one,
two, and three subspaces up to an isomorphism. A list of nonisomorphic
transitive systems of $n$ subspaces will be called complete if, for
any transitive system $S=(H;H_1,H_2,\ldots,H_n)$ of $n$ subspaces of
the space $H$, there is in the list  a system isomorphic to the system
$S$.

\begin{pr}\label{u2}
  If a system $S=(H;H_1)$ of a single subspace $H_1$ of the space $H$
  is transitive, then it is isomorphic to one of the following
  systems:
  $$
  S_1=(\mathbb C;0),\quad S_2=(\mathbb C;\mathbb C).
  $$
\end{pr}

\begin{proof}
  Let $\dim H>1$ and $H_1$ be an arbitrary proper subspace of the
  space $H$. Then the algebra $\End(S)$ corresponding to the system
  $S=(H;H_1)$ contains a nontrivial idempotent, for example, the
  operator of orthogonal projection onto $H^\bot_1$, and,
  consequently, the algebra is trivial. In the case where $\dim H>1$ and
  $H_1$ is a trivial subspace of the space $H$, the algebra
  $\End(S)=B(H)$, that is, it coincides with the set of linear bounded
  operators from $H$ into $H$.
\end{proof}

To construct lists of transitive systems of two and three subspaces,
we use the description of the algebra $\End(S)$ for the system
$S=(U;K_1,K_2,K_3)$ of $3$ subspaces $K_1$, $K_2$, $K_3$ of a finite
dimensional linear space $U$~\cite{B}. Let $L$ be an arbitrary
subspace complementary to the subspace $K_1+K_2+K_3$ in the space $U$,
that is, 
$$
(K_1+K_2+K_3)\dot{+}L=U,
$$
where $\dot{+}$ is the direct sum of vector spaces.

Denote $P=K_1\cap K_2\cap K_3$. Let Пусть $M_1$,$M_2$,$M_3$
be arbitrary subspaces complementary to the subspaces
$K_1\cap(K_2+K_3)$, $K_2\cap(K_1+K_3)$, $K_3\cap(K_1+K_2)$ in
$K_1$, $K_2$, $K_3$, correspondingly, that is,
$$
\begin{array}{l}
K_1\cap(K_2+K_3)\dot{+}M_1=K_1,\\
K_2\cap(K_1+K_3)\dot{+}M_2=K_2,\\
K_3\cap(K_1+K_2)\dot{+}M_3=K_3.
\end{array}
$$

Denote by $N_1$,$N_2$,$N_3$ arbitrary complementary subspaces to the
subspace $P$ in $K_2\cap K_3$, $K_1\cap
K_3$, $K_1\cap K_2$, correspondingly, that is,
$$
\begin{array}{l}
P\dot{+}N_1=K_2\cap K_3,\\
P\dot{+}N_2=K_1\cap K_3,\\
P\dot{+}N_3=K_1\cap K_2.
\end{array}
$$

Let now $Q_3$ be an arbitrary subspace complementary to the subspace
$K_3\cap K_1+K_3\cap K_1$ in the subspace $K_3\cap(K_1+K_2)$. An
arbitrary element $x_3$ of the subspace $Q_3$ is uniquely decomposed
into the sum $x_3=x_1+x_2$, where $x_1\in K_1$ and $x_2\in K_2$ are
such that if $x_3$ runs over a basis of $Q_3$, $x_1$ runs over a system
of linearly independent vectors the linear span of which makes a
subspace complementary to the subspace $K_1\cap K_2+K_1\cap K_3$ in
the space $K_1\cap(K_2+K_3)$, and $x_2$ runs over a system of linearly
independent vectors that span a subspace complementary to the subspace
$K_2\cap K_1+K_2\cap K_3$ in the subspace $K_2\cap(K_1+K_3)$. Denote
these complementary  subspaces by $Q_1$ and $Q_2$,
correspondingly. Thus,
$$
\begin{array}{l}
(K_1\cap K_2+K_1\cap K_3)\dot{+}Q_1=K_1\cap(K_2+K_3),\\
(K_2\cap K_1+K_2\cap K_3)\dot{+}Q_2=K_2\cap(K_1+K_3),\\
(K_3\cap K_1+K_3\cap K_2)\dot{+}Q_3=K_3\cap(K_1+K_2),
\end{array}
$$
and $\dim Q_1=\dim Q_2=\dim Q_3$. For the space $U$ and the subspaces
$K_1$, $K_2$, $K_3$, we have
\begin{equation}\label{e1}
\begin{array}{l}
U=L\dot{+}M_1\dot{+}M_2\dot{+}M_3\dot{+}Q_1\dot{+}Q_2\dot{+}N_1\dot{+}N_2\dot{+}N_3\dot{+}P,\\
K_1=M_1\dot{+}N_2\dot{+}N_3\dot{+}Q_1\dot{+}P,\\
K_2=M_2\dot{+}N_1\dot{+}N_3\dot{+}Q_2\dot{+}P,\\
K_3=M_3\dot{+}N_1\dot{+}N_2\dot{+}Q_3\dot{+}P.
\end{array}
\end{equation}
Let now $\ell$, $m_i$, $q$, $n_i$, $p$, $u$ be dimensions of $L$,
$M_i$, $Q_i$, $N_i$, $P$, and $U$, correspondingly. Then the dimension
of the algebra $\End(S)$ that corresponds to the system
$S=(U;K_1,K_2,K_3)$, considered as a linear space, can be calculated
by the formula
\begin{equation}\label{e2}
\begin{array}{l}
\dim\,\End(S)=\ell
u+q^2+q\sum\limits_{i=1}^3(m_i+n_i)+\sum\limits_{i=1}^3(m_i^2+n_i^2)+\\
+\sum\limits_{\substack{ i\neq j \\ i,j=1}}^3m_in_j+p^2.
\end{array}
\end{equation}

\begin{pr}\label{u3}
  If a system $S=(H;H_1,H_2)$ of two subspaces of a space $H$ is
  transitive, then it is isomorphic to one of the following system:
  $$
  \begin{array}{cc}
    S_1=(\mathbb C;0,0),&  S_3=(\mathbb C;0,\mathbb C),\\
    S_2=(\mathbb C;\mathbb C,0),& S_4=(\mathbb C;\mathbb C,\mathbb C).
  \end{array}
  $$
\end{pr}

\begin{proof}
  To make an analysis of a system of two subspaces in the case of a
  finite dimensional linear space, set  $U=H$,
  $K_1=H_1$, $K_1=H_1$, $K_3=0$ in identities~(\ref{e1}). We get
  $$
  \begin{array}{l}
    H=L\dot{+}M_1\dot{+}M_2\dot{+}N_3,\\
    H_1=M_1\dot{+}N_3,\\ H_2=M_2\dot{+}N_3.
  \end{array}
  $$
  The formula for the dimension of the algebra $\End(S)$, for $K_3=0$, becomes
  $$
  \dim\, \End(S)=\ell u+m_1^2+m_2^2+n_3^2.
  $$
  Since the system $S=(H;H_1,H_2)$ is transitive, it follows that
  $\dim\End(S)=~1$ and, correspondingly, $\ell u+m_1^2+m_2^2+n_3^2=1$.
  This identity can hold only in the following four cases:
  \begin{itemize}
  \item[1)]$\ell u=1$. Hence,  $\dim L=1$, $H=L$, $H_1=0$,
    $H_2=0$ and, consequently, $S\cong S_1$.
  \item[2)]$m_1^2=1$. Hence, $\dim M_1=1$, $H=M_1$,
    $H_1=M_1$, $H_2=0$ and, consequently, $S\cong S_2$.
  \item[3)]$m_2^2=1$. Hence, $\dim M_2=1$, $H=M_2$,
    $H_1=0$, $H_2=M_2$ and, consequently, $S\cong S_3$.
  \item[4)]$n_3^2=1$. Hence, $\dim N_3=1$, $H=N_3$,
    $H_1=N_3$, $H_2=N_3$ and, consequently, $S\cong S_4$.
  \end{itemize}
  It follows from Proposition~\ref{u1} and~\cite{OS} that if a pair of
  orthogonal projections on an infinite dimensional Hilbert space is
  reducible, then there do not exist transitive systems of two
  subspaces in an infinite dimensional Hilbert space. We remark that
  this fact can also be obtained from decomposability of a system of two
  subspaces in an infinite dimensional Hilbert space~\cite{EW}.
\end{proof}

\begin{pr}\label{u4}
  If a system $S=(U;K_1,K_2,K_3)$ of three subspaces of a finite
  dimensional linear space $U$ is transitive, then it is isomorphic to
  one of the following systems:
  $$
  \begin{array}{c}
    \begin{array}{cc}
      S_1=(\mathbb C;0,0,0),&  S_5=(\mathbb C;0,\mathbb C,\mathbb C),\\
      S_2=(\mathbb C;\mathbb C,0,0),& S_6=(\mathbb C;\mathbb C,0,\mathbb
      C),\\
      S_3=(\mathbb C;0,\mathbb C,0),& S_7=(\mathbb C;\mathbb C,\mathbb
      C,0),\\
      S_4=(\mathbb C;0,0,\mathbb C),& S_8=(\mathbb C;\mathbb C,\mathbb
      C,\mathbb C),
    \end{array}\\
    S_9=(\mathbb C^2;\mathbb C(1,0),\mathbb C(0,1),\mathbb C(1,1)).
  \end{array}
  $$
\end{pr}

\begin{proof}
  Since the system $S=(U;K_1,K_2,K_3)$ is transitive, it follows that
  $\dim\End(S)=~1$ and, correspondingly,
  $$
    \ell
    u+q^2+q\sum\limits_{i=1}^3(m_i+n_i)+\sum\limits_{i=1}^3(m_i^2+n_i^2)
    +\sum\limits_{\substack{i\neq j \\ i,j=1}}^3m_in_j+p^2=1.
  $$
  The last identity can hold only in one of the following nine
  cases:
  \begin{itemize}
  \item[1)]$\ell u=1$. Hence, $\dim L=1$, $U=L$, $K_1=0$,
    $K_2=0$, $K_3=0$. Thus  $S\cong S_1$.
  \item[2)]$m_1^2=1$. Hence, $\dim M_1=1$, $U=M_1$,
    $K_1=M_1$, $K_2=0$, $K_3=0$ and thus $S\cong S_2$.
  \item[3)]$m_2^2=1$. Hence, $\dim M_2=1$, $U=M_2$,
    $K_1=0$, $K_2=M_2$, $K_3=0$, and thus $S\cong S_3$.
  \item[4)]$m_3^2=1$. Hence, $\dim M_3=1$, $U=M_3$,
    $K_1=0$, $K_2=0$, $K_3=M_3$, and thus $S\cong S_4$.
  \item[5)]$n_1^2=1$. Hence, $\dim N_1=1$, $U=N_1$,
    $K_1=0$, $K_2=N_1$, $K_3=N_1$, and thus $S\cong S_5$.
  \item[6)]$n_2^2=1$. Hence, $\dim N_2=1$, $U=N_2$,
    $K_1=N_2$, $K_2=0$, $K_3=N_2$, and thus $S\cong S_6$.
  \item[7)]$n_3^2=1$. Hence, $\dim N_3=1$, $U=N_3$,
    $K_1=N_3$, $K_2=N_3$, $K_3=0$, and thus $S\cong S_7$.
  \item[8)]$p^2=1$. Hence, $\dim P=1$, $U=P$,
    $K_1=P$, $K_2=P$, $K_3=P$, and thus $S\cong S_8$.
  \item[9)]$q^2=1$. Hence, $\dim Q_1=\dim Q_2=1$, $U=Q_1\dot{+}Q_2$,
    $K_1=Q_1$, $K_2=Q_2$, $K_3=Q_3$, and thus $S\cong S_9$.
  \end{itemize}
\end{proof}

We recall that the problem of even proving existence of transitive
triples of subspaces of an infinite dimensional space is an open
problem (see~\cite{Ha2}).

\subsection{Transitive systems of four subspaces}

Following~\cite{GP} let us introduce the notion of a  defect of a
system $S=(U;K_1,K_2,K_3,K_4)$ of four subspaces of a finite
dimensional linear space $U$.

\begin{defn}\label{d10}
  Let $S=(U;K_1,K_2,K_3,K_4)$ be a system of four subspaces of a
  finite dimensional linear space $U$. By a defect of the system $S$,
  we will call the number defined by
  $$
  \rho(S)=\sum\limits_{i=1}^4\dim K_i-2\dim U.
  $$
\end{defn}

S.~Brenner in~\cite{B} gave a description of a complete list of four
distinct proper subspaces up to a rearrangement of the subspaces, and
systems that have a nonnegative defect were written down explicitly.
An explicit form for systems of four proper subspaces, with a negative
defect, is given in this section by passing to orthogonal systems and
choosing suitable isomorphic systems. We adopt the following notations
used in~\cite{B}:
\begin{itemize}
\item[] $\bf 1$ is the $r\times r$ identity matrix;
\item[] $\bf 0$ is the $r\times r$ zero matrix;
\item[] $\bf J$ is the $r\times r$ Jordan cell with zero on the diagonal;
\item[] $\xi$ is the column of $r$ zeros;
\item[] $\eta$ is the row of $r$ zeros;
\item[] $b$ is the column of the first $(r-1)$ zeros and $1$ as the last
  element;
\item[] $d$ is the row with the first element equal $1$ and other
  $r-1$ zeros.
\end{itemize}

The subspace $K_i$ in the list is given by a matrix $\EuScript
K_i$. Here the subspace $K_i$ is set to be the linear span of rows of
the matrix $\EuScript K_i$. Introduce two more notations, ---
$B(u,\rho)$ denotes the system $B=(U;K_1,K_2,K_3,K_4)$ of four
subspaces of the space $U$ of dimension $u$ with defect $\rho$, and
$B(u,\rho;\lambda)$ denotes the system $B=(U;K_1,K_2,K_3,K_4)$ of four
subspaces of the spaces $U$ of dimension $u$, with defect $\rho$,
which depend on a parameter $\lambda$.

The following is a complete list of distinct proper subspaces, up to a
rearrangement:

(1) $B(2,0;\lambda)$, $\lambda\in \mathbb C, \lambda\neq0,1$,
$$
\EuScript K_1=\begin{pmatrix} 1& 0 \end{pmatrix},\quad \EuScript
K_2=\begin{pmatrix} 0& 1 \end{pmatrix},\quad\EuScript
K_3=\begin{pmatrix} 1&1 \end{pmatrix},\quad \EuScript
K_4=\begin{pmatrix} 1& \lambda \end{pmatrix}.
$$

(2) $B(2r,1)$, $r=2,3,\ldots$,
$$
\EuScript K_1=\begin{pmatrix} \bf{1}& \bf{0} \end{pmatrix},\quad
\EuScript K_2=\begin{pmatrix} \bf{0}& \bf{1}
\end{pmatrix},\quad\EuScript K_3=\begin{pmatrix} \bf{1}&\bf{1}
\end{pmatrix},\quad \EuScript K_4=\begin{pmatrix} \bf{1}&
\bf{J}\\\eta & d \end{pmatrix}.
$$

(3) $B(2r+2,-1)$, $r=1,2,\ldots$,
$$
\EuScript K_1=\begin{pmatrix} \bf{1}& \bf{0} & \xi & \xi\\ \eta &
d &0&0 \end{pmatrix},\quad \EuScript K_2=\begin{pmatrix} \bf{0} &
\bf{J} & b & \xi\\ \eta & \eta & 0 & 1\end{pmatrix},
$$
$$
\EuScript K_3=\begin{pmatrix} \bf{1} & \bf{J} & b & \xi\\ \eta & d
& 0& 1\end{pmatrix},\quad \EuScript K_4=\begin{pmatrix} \bf{1}&\xi
& \xi & \bf{1} \end{pmatrix}.
$$

(4a) $B(3,1)$,
$$
\EuScript K_1=\begin{pmatrix}1&0&0\\0&1&0\end{pmatrix},\quad
\EuScript K_2=\begin{pmatrix}1&0&0\\0&0&1\end{pmatrix},
$$
$$
\EuScript K_3=\begin{pmatrix}0&1&0\\0&0&1\end{pmatrix},\quad
\EuScript K_4=\begin{pmatrix}1&1&1\end{pmatrix}.
$$

(4b) $B(2r+3,1)$, $r=1,2,\ldots$,
$$
\EuScript K_1=\begin{pmatrix}\bf{1}&\bf{0}&\xi&\xi&\xi\\
\eta&\eta&1&0&0\\\eta&\eta&0&1&0 \end{pmatrix},\quad \EuScript
K_2=\begin{pmatrix}\bf{0}&\bf{1}&\xi&\xi&\xi\\
\eta&\eta&1&0&0\\\eta&\eta&0&0&1\end{pmatrix},
$$
$$
\EuScript K_3=\begin{pmatrix}\bf{1}&\bf{1}&\xi&\xi&\xi\\
\eta&\eta&0&1&0\\\eta&\eta&0&0&1\end{pmatrix},\quad \EuScript
K_4=\begin{pmatrix}\bf{1}&\bf{J}&b&\xi&b\\ \eta & d &
0&1&0\end{pmatrix}.
$$

(5a) $B(3,-1)$,
$$
\EuScript K_1=\begin{pmatrix}0&1&0\end{pmatrix},\quad \EuScript
K_2=\begin{pmatrix}0&0&1\end{pmatrix},
$$
$$
\EuScript K_3=\begin{pmatrix}1&0&0\end{pmatrix},\quad \EuScript
K_4=\begin{pmatrix}0&1&1\\1&0&1\end{pmatrix}.
$$

(5b) $B(2r+3,-1)$, $r=1,2,\ldots$,
$$
\EuScript K_1=\begin{pmatrix}\bf{1}&\bf{0}&\xi&\xi&\xi\\
\eta&\eta&0&1&0\end{pmatrix},\quad \EuScript
K_2=\begin{pmatrix}\bf{0}&\bf{1}&\xi&\xi&\xi\\
\eta&\eta&0&0&1\end{pmatrix},
$$
$$
\EuScript K_3=\begin{pmatrix}\bf{1}&\bf{1}&\xi&\xi&\xi\\
\eta&\eta&1&0&0\end{pmatrix},\quad \EuScript
K_4=\begin{pmatrix}\bf{1}&\bf{J}&b&\xi&\xi\\ \eta & d & 0&1&0
\\ \eta&\eta&1&0&1\end{pmatrix}.
$$

(6a) $B(3,2)$,
$$
\EuScript K_1=\begin{pmatrix}1&0&0\\ 0&1&0\end{pmatrix}, \quad
\EuScript K_2=\begin{pmatrix}1&0&0 \\ 0&0&1 \end{pmatrix},
$$
$$
\EuScript K_3=\begin{pmatrix}0&1&0 \\ 0&0&1 \end{pmatrix},\quad
\EuScript K_4=\begin{pmatrix}1&0&1\\1&1&0\end{pmatrix}.
$$

(6b) $B(5,2)$,
$$
\EuScript K_1=\begin{pmatrix}1&0&0&0&0\\0& 0&1&0&0
\\ 0&0&0&1&0\end{pmatrix}, \quad \EuScript K_2=\begin{pmatrix}0&1&0&0&0 \\
0&0&1&0&0\\ 0&0&0&0&1 \end{pmatrix},
$$
$$
\EuScript K_3=\begin{pmatrix}1&1&0&0&0 \\ 0&0&0&1&0 \\ 0&0&0&0&1
\end{pmatrix},\quad \EuScript
K_4=\begin{pmatrix}1&0&1&0&0 \\ 0&1&0&0&0 \\
0&0&1&1&1 \end{pmatrix}.
$$

(6c) $B(2r+3,2)$, $r=2,3,\ldots$,
$$
\EuScript K_1=\begin{pmatrix}\bf{1}&\bf{0}&\xi&\xi&\xi\\\eta&
\eta&1&0&0
\\ \eta&\eta&0&1&0\end{pmatrix}, \quad \EuScript K_2=\begin{pmatrix}\bf{0}&\bf{1}&\xi&\xi&\xi \\
\eta&\eta&1&0&0\\ \eta&\eta&0&0&1 \end{pmatrix},
$$
$$
\EuScript K_3=\begin{pmatrix}\bf{1}&\bf{1}&\xi&\xi&\xi \\ \eta&\eta&0&1&0 \\
\eta&\eta&0&0&1
\end{pmatrix},\quad \EuScript
K_4=\begin{pmatrix}\bf{1}&{\bf J}^2&{\bf J}b&\xi&({\bf J}+{\bf 1})b \\ \eta&d&0&0&0 \\
\eta&d{\bf J}&0&1&0 \end{pmatrix}.
$$

(7a) $B(3,-2)$,
$$
\EuScript K_1=\begin{pmatrix}0&1&0\end{pmatrix}, \quad \EuScript
K_2=\begin{pmatrix}0&0&1 \end{pmatrix},\quad \EuScript
K_3=\begin{pmatrix}1&0&0 \end{pmatrix},\quad \EuScript
K_4=\begin{pmatrix}1&1&1\end{pmatrix}.
$$

(7b) $B(5,-2)$,
$$
\EuScript K_1=\begin{pmatrix}1&0&0&0&0
\\ 0&0&0&1&0\end{pmatrix}, \quad \EuScript K_2=\begin{pmatrix}0&1&0&0&0 \\
 0&0&0&0&1 \end{pmatrix},
$$
$$
\EuScript K_3=\begin{pmatrix}1&1&0&0&0 \\ 0&0&1&0&0
\end{pmatrix},\quad \EuScript
K_4=\begin{pmatrix}1&0&1&1&0 \\  0&0&0&1&1 \end{pmatrix}.
$$

(7c) $B(2r+5,-2)$, $r=1,2,\ldots$,
$$
\EuScript K_1=\begin{pmatrix}{\bf 1}&{\bf 0}&\xi&\xi&\xi&\xi&\xi
\\\eta& d&0&0&0&0&0
\\ \eta&\eta&0&0&0&1&0\end{pmatrix}, \quad \EuScript K_2=
 \begin{pmatrix}{\bf 0}&{\bf J}&b&\xi&\xi&\xi&\xi \\
\eta&\eta&0&1&0&0&0\\ \eta&\eta0&0&0&0&1 \end{pmatrix},
$$
$$
\EuScript K_3=\begin{pmatrix}{\bf 1}&{\bf J}&b&\xi&\xi&\xi&\xi \\ \eta& d & 0&1&0&0&0 \\
\eta&\eta&0&0&1&0&0
\end{pmatrix},\quad \EuScript
K_4=\begin{pmatrix}{\bf 1}&{\bf J}^3&{\bf J}^2b&{\bf J}b&b&\xi&\xi \\ b^T&d&0&0&0&0&1\\
\eta&d{\bf J}^2&0&0&0&1&0 \end{pmatrix}.
$$

\begin{thm}[S. Brenner]\label{t1}
  If a system $S=(U;K_1,K_2,K_3,K_4)$ of four distinct proper
  subspaces of a finite dimensional linear space $U$ is transitive,
  then it is isomorphic, up to a rearrangement of the subspaces, to
  one of the following system:
  $$
  \begin{array}{c}
    B(2,0;\lambda), \quad\lambda\in \mathbb C, \lambda\neq0,1,\\
    B(u,\pm1), \quad u=3,4,5,\ldots,\\
    B(u,\pm2), \quad u=3,5,7,\ldots.
  \end{array}
  $$
\end{thm}

\section{The algebra $\EuScript P_{n,com}$ and its
  $*$-representations}

\subsection{Irreducible $*$-representations of the algebra $\EuScript
  P_{n,com}$} 

For $n\in\mathbb N$, denote by $\Sigma_n$ the set of $\alpha\in\mathbb
R_+$ such that there exists at least one $*$-representation of the
\hbox{$*$-algebra} $\EuScript P_{n,\alpha}=\mathbb C<p_1,p_2,\ldots,
p_n|p^2_k=p^*_k=p_k,\sum^n_{k=1}p_k=\alpha e>$, that is, the set of
all real parameters $\alpha$ for which there exist $n$ orthogonal
projections $P_1$, $P_2$, $\ldots$, $P_n$ on a Hilbert space $H$
satisfying the relation $\sum^n_{k=1}P_k=\alpha I_H$. Introduce an
algebra, $\EuScript P_{n,com}=\mathbb C<p_1,p_2,\ldots, p_n|p^2_k=p^*_k=p_k,
[\sum_{k=1}^np_k,p_i]=0,\,\forall i=\overline{1,n}>$. All irreducible
$*$-representations of $\EuScript P_{n,com}$ is a union over all
$\alpha\in\Sigma_n$ of irreducible $*$-representations of $\EuScript
P_{n,\alpha}$.

A description of the set $\Sigma_n$ for all $n\in\mathbb N$ was
obtained by S.~A.~Kruglyak, V.~I.~Rabanovich, and
Yu.~S.~Samo\v\i{}lenko in~\cite{KRS}, and is given by
$$
\begin{array}{c}
\displaystyle
\Sigma_1=\{0,1\},\quad\Sigma_2=\{0,1,2\},\quad\Sigma_3=\{0,1,
{\textstyle\frac32},2,3\},
\\[4mm]
\displaystyle
\Sigma_n=\{\Lambda_n^0,\Lambda_n^1,
\left[{\textstyle\frac{n-\sqrt{n^2-4n}}2,\frac{n+\sqrt{n^2-4n}}2}\right],
n-\Lambda_n^1,n-\Lambda_n^0\}, \,\mbox{ при}\, n\geq4,
\\[4mm]
\displaystyle
\Lambda_n^0=\Biggl\{
 0,1+{\textstyle \frac1{n-1}}, 1+{\textstyle
\frac1{(n-2)-\frac1{n-1}}},\ldots, 1+{\textstyle
\frac1{(n-2)-\frac1{(n-2)-\frac1{\ddots-\frac1{n-1}}}}},\ldots
 \Biggr \},
\\[4mm]
\displaystyle
\Lambda_n^1=\Biggl\{
 0,1+{\textstyle \frac1{n-2}}, 1+{\textstyle
\frac1{(n-2)-\frac1{n-2}}},\ldots, 1+{\textstyle
\frac1{(n-2)-\frac1{(n-2)-\frac1{\ddots-\frac1{n-2}}}}},\ldots
 \Biggr \}.
\end{array}
$$

\subsection{Irreducible $*$-representations of the algebras $\EuScript
P_{1,com}$, $\EuScript P_{2,com}$, $\EuScript P_{3,com}$} 

Let us give a list of irreducible $*$-representations of the algebra
$\EuScript P_{1,com}$. By~\cite{KRS}, we have $\Sigma_1=\{0,1\}$.

For $\alpha=0$, the only irreducible representation of the algebra
$\EuScript P_{1,0}$, up to equivalence, is the representation $P_1=0$
on the space $H=\mathbb C$. For $\alpha=1$, the unique up to
equivalence irreducible representation of the algebra $\EuScript
P_{1,1}$ is the representation $P_1=\mathbb C$ on the space $H=\mathbb C$.

For the algebra $\EuScript P_{2,com}$, we have
$\Sigma_2=\{0,1,2\}$~\cite{KRS}.

If $\alpha=0$, there is a unique up to equivalence irreducible
representation of the algebra $\EuScript P_{2,0}$ given by $P_1=0$,
$P_2=0$ on the space $H=\mathbb C$. If $\alpha=1$, there are two
irreducible representations of the algebra $\EuScript P_{2,1}$, up to
equivalence. The first one is given by $P_1=I$, $P_2=0$ on the space
$H=\mathbb C$, and the second one by $P_1=0$, $P_2=I$ on the space
$H=\mathbb C$.  In the case where $\alpha=2$, the only representation
of the algebra $\EuScript P_{2,2}$, up to equivalence, is the
representation $P_1=I$, $P_2=I$ on the space $H=\mathbb C$.

Now we give irreducible $*$-representations of the algebra $\EuScript
P_{3,com}$. We have $\Sigma_3=\{0,1,
{\textstyle\frac32},2,3\}$.

If $\alpha=0$, there is a unique up to equivalence irreducible
representation of the algebra $\EuScript P_{3,0}$. It is given by
$P_1=0$, $P_2=0$, $P_3=0$ on the space $H=\mathbb
C$. If $\alpha=1$, there are three inequivalent irreducible
representations of the algebra $\EuScript P_{3,1}$. The first one is
$P_1=I$, $P_2=0$, $P_3=0$ on the space $H=\mathbb C$. The second one
is $P_1=0$, $P_2=I$, $P_3=0$ on $H=\mathbb C$. The third one is given
by $P_1=0$, $P_2=0$, $P_3=I$ on the space $H=\mathbb C$. If
$\alpha=3/2$, there is a unique up to equivalence irreducible
representation of the algebra $\EuScript P_{3,3/2}$,
$$
P_1=\begin{pmatrix}1&0\\0&0\end{pmatrix},\quad
P_2=\begin{pmatrix}1/4&\sqrt{3}/4\\\sqrt{3}/4&3/4\end{pmatrix},\quad
P_3=\begin{pmatrix}1/4&-\sqrt{3}/4\\-\sqrt{3}/4&3/4\end{pmatrix},
$$
which act on the space $H=\mathbb C^2$. If $\alpha=2$, there are
three inequivalent irreducible representations of the algebra
$\EuScript P_{3,2}$. The first one is $P_1=0$, $P_2=I$, $P_3=I$ on
$H=\mathbb C$, the second one is $P_1=I$, $P_2=0$, $P_3=I$ on
$H=\mathbb C$, and the third one is $P_1=I$, $P_2=I$, $P_3=0$ on
$H=\mathbb C$. For $\alpha=3$, the unique up to equivalence
irreducible representation of the algebra $\EuScript P_{3,3}$ is
$P_1=I$, $P_2=I$, $P_3=I$ on $H=\mathbb C$.

\subsection{Irreducible $*$-representations of the algebra $\EuScript
P_{4,com}$} 

We use the following notations:
$$
\begin{array}{l}
\displaystyle
A_{\ell,m}=\frac1m\begin{pmatrix}m-\ell&-\sqrt{\ell(m-l)}\\
-\sqrt{\ell(m-l)}&\ell\end{pmatrix},\\[4mm]
\displaystyle
B_{\ell,m}=\frac1m\begin{pmatrix}m-\ell&\sqrt{\ell(m-l)}\\
\sqrt{\ell(m-l)}&\ell\end{pmatrix},\\[4mm]
\displaystyle
C_{\ell,m}=I-A_{\ell,m}
=\frac1m\begin{pmatrix}\ell&\sqrt{\ell(m-l)}\\
\sqrt{\ell(m-l)}&m-\ell\end{pmatrix},\\[4mm]
\displaystyle
D_{\ell,m}=I-B_{\ell,m}
=\frac1m\begin{pmatrix}\ell&-\sqrt{\ell(m-l)}\\
-\sqrt{\ell(m-l)}&m-\ell\end{pmatrix}.
\end{array}
$$

Let us consider a part of the unit sphere $\Omega\subset\mathbb R^3$,
given by $\Omega=\{(a,b,c)\in\mathbb R| a^2+b^2+c^2=1, a>0, b>0,
c\in(-1,1)\, \mbox{или}\, a=0, b^2+c^2=1, b>0, c>0 \, \mbox{или}\,
b=0, a^2+c^2=1, b>0, c>0\}$.

\begin{figure}[h]
  \includegraphics[height=6cm] {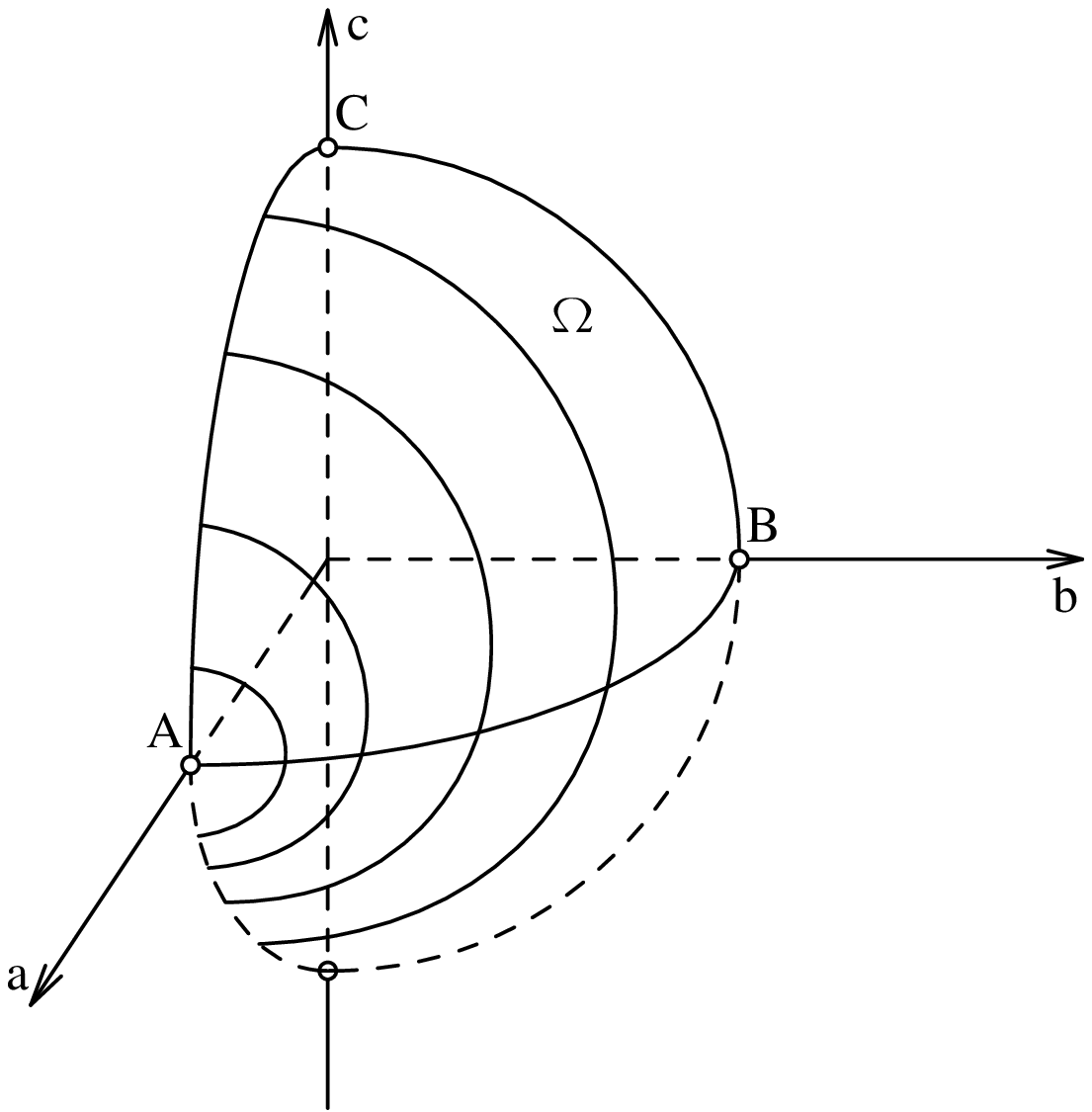}
  \caption{}\label{fig1}
\end{figure}

Since all irreducible $*$-representations of the algebra $\EuScript
P_{4,com}$ are finite dimensional,  denote the space of
representations by $U$. Also denote by $S(u,\rho)$ the system $S=(U;\Imp
P_1,\Imp P_2,\Imp P_3,\Imp P_4)$ of four subspaces of the space $U$ of
dimension $u$ with defect $\rho$, which is generated by the
representation $P_1$, $P_2$, $P_3$, $P_4$ on the space $U$, and by
$S(u,\rho;a,b,c)$ the systems $S=(U;\Imp P_1,\Imp P_2,\Imp P_3,\Imp
P_4)$ of four subspaces of the space $U$ of dimension $u$ with
defect $\rho$, which are generated by the representation $P_1$,
$P_2$,$P_3$, $P_4$ on $U$ and depend on the  parameters
$a$,$b$,$c$. Using the results of~\cite{KRS,OS}, we write a list
of systems of four distinct proper subspaces, given up to a
rearrangement of the subspaces, which are generated by irreducible
inequivalent representations of the algebra $\EuScript
P_{4,\alpha}$:

(1) $S(2,0;a,d,c)$, $(a,b,c)\in\Omega$,
$$
P_1=\frac12\begin{pmatrix}1+a&-b-ic\\-b+ic&1-a\end{pmatrix},\,
P_3=\frac12\begin{pmatrix}1-a&-b+ic\\-b-ic&1+a\end{pmatrix},
$$
$$
P_2=\frac12\begin{pmatrix}1-a&b-ic\\b+ic&1+a\end{pmatrix},\qquad
P_4=\frac12\begin{pmatrix}1+a&b+ic\\b-ic&1-a\end{pmatrix}.
$$

(2) $S(2r,1)$, $r=2,3,\ldots$,
$$
\begin{array}{l}
P_1=A_{2r-1,4r}\oplus A_{2r-3,4r}\oplus\ldots\oplus A_{1,4r},\\
P_2=B_{2r-1,4r}\oplus B_{2r-3,4r}\oplus\ldots\oplus B_{1,4r},\\
U=\underbrace{\mathbb C^2\oplus\ldots\oplus\mathbb C^2}_r;
\end{array}
$$
$$
\begin{array}{l}
P_3=0\oplus B_{2r-2,4r}\oplus B_{2r-4,4r}\ldots\oplus B_{2,4r}\oplus1,\\
P_4=1\oplus A_{2r-2,4r}\oplus A_{2r-4,4r}\ldots\oplus A_{2,4r}\oplus1,\\
U=\mathbb C\oplus\underbrace{\mathbb C^2\oplus\ldots\oplus\mathbb
C^2}_{r-1}\oplus\mathbb C ;
\end{array}
$$

(3) $S(2r,-1)$, $r=2,3,\ldots$,
$$
\begin{array}{l}
P_1=C_{2r-1,4r}\oplus C_{2r-3,4r}\oplus\ldots\oplus C_{1,4r},\\
P_2=D_{2r-1,4r}\oplus D_{2r-3,4r}\oplus\ldots\oplus D_{1,4r},\\
U=\underbrace{\mathbb C^2\oplus\ldots\oplus\mathbb C^2}_r;
\end{array}
$$
$$
\begin{array}{l}
P_3=1\oplus D_{2r-2,4r}\oplus D_{2r-4,4r}\ldots\oplus D_{2,4r}\oplus0,\\
P_4=0\oplus C_{2r-2,4r}\oplus C_{2r-4,4r}\ldots\oplus C_{2,4r}\oplus0,\\
U=\mathbb C\oplus\underbrace{\mathbb C^2\oplus\ldots\oplus\mathbb
C^2}_{r-1}\oplus\mathbb C ;
\end{array}
$$

(4) $S(2r+1,1)$, $r=1,2,\ldots$,
$$
\begin{array}{l}
P_1=A_{2r,4r+2}\oplus A_{2r-2,4r+2}\oplus\ldots\oplus A_{2,4r+2}\oplus1,\\
P_2=B_{2r,4r+2}\oplus B_{2r-2,4r+2}\oplus\ldots\oplus B_{2,4r+2}\oplus1,\\
U=\underbrace{\mathbb C^2\oplus\ldots\oplus\mathbb
C^2}_r\oplus\mathbb C;
\end{array}
$$
$$
\begin{array}{l}
P_3=1\oplus B_{2r-1,4r+2}\oplus B_{2r-3,4r+2}\ldots\oplus B_{1,4r+2},\\
P_4=0\oplus A_{2r-1,4r+2}\oplus A_{2r-3,4r+2}\ldots\oplus A_{1,4r+2},\\
U=\mathbb C\oplus\underbrace{\mathbb C^2\oplus\ldots\oplus\mathbb
C^2}_r.
\end{array}
$$

(5) $S(2r+1,-1)$, $r=1,2,\ldots$,
$$
\begin{array}{l}
P_1=C_{2r,4r+2}\oplus C_{2r-2,4r+2}\oplus\ldots\oplus C_{2,4r+2}\oplus0,\\
P_2=D_{2r,4r+2}\oplus D_{2r-2,4r+2}\oplus\ldots\oplus D_{2,4r+2}\oplus0,\\
U=\underbrace{\mathbb C^2\oplus\ldots\oplus\mathbb
C^2}_r\oplus\mathbb C;
\end{array}
$$
$$
\begin{array}{l}
P_3=0\oplus D_{2r-1,4r+2}\oplus D_{2r-3,4r+2}\ldots\oplus D_{1,4r+2},\\
P_4=1\oplus C_{2r-1,4r+2}\oplus C_{2r-3,4r+2}\ldots\oplus C_{1,4r+2},\\
U=\mathbb C\oplus\underbrace{\mathbb C^2\oplus\ldots\oplus\mathbb
C^2}_r.
\end{array}
$$

(6) $S(2r+1,2)$, $r=1,2,\ldots$,
$$
\begin{array}{l}
P_1=1\oplus A_{2r-1,2r+1}\oplus A_{2r-3,2r+1}\oplus\ldots\oplus A_{1,2r+1},\\
P_2=1\oplus B_{2r-1,2r+1}\oplus B_{2r-3,2r+1}\oplus\ldots\oplus B_{1,2r+1},\\
U=\mathbb C\oplus\underbrace{\mathbb C^2\oplus\ldots\oplus\mathbb
C^2}_r;
\end{array}
$$
$$
\begin{array}{l}
P_3=B_{2r,2r+1}\oplus B_{2r-2,2r+1}\ldots\oplus B_{2,2r+1}\oplus1,\\
P_4=A_{2r,2r+1}\oplus A_{2r-2,2r+1}\ldots\oplus A_{2,2r+1}\oplus1,\\
U=\underbrace{\mathbb C^2\oplus\ldots\oplus\mathbb
C^2}_r\oplus\mathbb C.
\end{array}
$$

(7) $S(2r+1,-2)$, $r=1,2,\ldots$,
$$
\begin{array}{l}
P_1=0\oplus C_{2r-1,2r+1}\oplus C_{2r-3,2r+1}\oplus\ldots\oplus C_{1,2r+1},\\
P_2=0\oplus D_{2r-1,2r+1}\oplus D_{2r-3,2r+1}\oplus\ldots\oplus D_{1,2r+1},\\
U=\mathbb C\oplus\underbrace{\mathbb C^2\oplus\ldots\oplus\mathbb
C^2}_r;
\end{array}
$$
$$
\begin{array}{l}
P_3=D_{2r,2r+1}\oplus D_{2r-2,2r+1}\ldots\oplus D_{2,2r+1}\oplus0,\\
P_4=C_{2r,2r+1}\oplus C_{2r-2,2r+1}\ldots\oplus C_{2,2r+1}\oplus0,\\
U=\underbrace{\mathbb C^2\oplus\ldots\oplus\mathbb
C^2}_r\oplus\mathbb C.
\end{array}
$$

Hence, irreducible inequivalent representations, $\Rep\EuScript
P_{4,\alpha}$, give rise to the following list of systems of four
distinct proper subspaces:
\begin{equation}\label{e3}
\begin{array}{c}
S(2,0;a,b,c), \, (a,b,c)\in \Omega,\\
S(u,\pm1), \quad u=3,4,5,\ldots,\\
S(u,\pm2), \quad u=3,5,7,\ldots.
\end{array}
\end{equation}

\section{Systems of subspaces generated by $\Rep\EuScript
P_{n,com}$, and transitive systems of $n$ subspaces}

\subsection{Transitive systems of subspaces generated by $\Rep\EuScript
P_{1,com}$, $\Rep\EuScript P_{2,com}$, $\Rep\EuScript
P_{3,com}$}

In this section we show that irreducible nonequivalent
$*$-representations of the $*$-algebras $\EuScript P_{1,com}$ and
$\EuScript P_{2,com}$ generate all nonisomorphic transitive systems of
one and two subspaces of an arbitrary Hilbert space. If $n=3$,
irreducible nonequivalent $*$-representations of the $*$-algebra
$\EuScript P_{3,com}$ give rise to all nonisomorphic transitive
systems of three subspaces of a finite dimensional linear space.

\begin{pr}\label{u5}
  Irreducible nonequivalent $*$-representations of $\EuScript
  P_{1,com}$ generate all transitive systems of one subspace of a
  Hilbert space.
\end{pr}

\begin{proof}
  Using Proposition~\ref{u2} we get a complete list of transitive
  systems of one subspaces as follows:
  $$
  S_1=(\mathbb C;0),\quad S_2=(\mathbb C;\mathbb C).
  $$
  By the results of Section~3, we have $\Sigma_1=\{0,1\}$.
  
  If $\alpha=0$, a unique up to equivalence irreducible representation
  of the algebra $\EuScript P_{1,0}$ is the representation $P_1=0$ on
  the space $H=\mathbb C$ and, consequently, a system of one subspace,
  induced by this representation, is isomorphic to $S_1$.
  
  If $\alpha=1$, there is only one, up to equivalence, irreducible
  representation of $\EuScript P_{1,1}$, $P_1=\mathbb C$, on the space
  $H=\mathbb C$, and so a system of one subspace, corresponding to
  this representation, is isomorphic to $S_2$.
\end{proof}

\begin{pr}\label{u6}
  Irreducible nonequivalent $*$-representations of $\EuScript
  P_{2,com}$ generate all transitive systems of two subspaces of a
  Hilbert space.
\end{pr}

\begin{proof}
  By Proposition~\ref{u3}, a complete list of transitive systems of
  two subspaces has the form
  $$
  \begin{array}{cc}
    S_1=(\mathbb C;0,0),&  S_3=(\mathbb C;0,\mathbb C),\\
    S_2=(\mathbb C;\mathbb C,0),& S_4=(\mathbb C;\mathbb C,\mathbb C).
  \end{array}
  $$
  By Section~3, $\Sigma_2=\{0,1,2\}$.
  
  For $\alpha=0$, the algebra $\EuScript P_{2,0}$ has, up to
  equivalence, a unique irreducible representation $P_1=0$,
  $P_2=0$ on the space $H=\mathbb C$ and, consequently, the system of
  subspaces generated by this representation is isomorphic to $S_1$.
  
  If $\alpha=1$, there are two inequivalent representations of
  $\EuScript P_{2,1}$. The first one is $P_1=I$, $P_2=0$ on the space
  $H=\mathbb C$. A system of two subspaces that corresponds to this
  representation is isomorphic to $S_2$. The second representation is
  given by $P_1=0$, $P_2=I$ on the space $H=\mathbb C$. The
  corresponding system of two subspaces is isomorphic to $S_3$.
  
  If $\alpha=2$, the only irreducible representation of the algebra
  $\EuScript P_{2,2}$ is $P_1=I$, $P_2=I$ on $H=\mathbb C$ and,
  consequently, the corresponding system of two subspaces is
  isomorphic to $S_4$.
\end{proof}

\begin{pr}\label{u7}
  Irreducible nonequivalent $*$-representations of $\EuScript
  P_{3,com}$ generate all transitive systems of three subspaces of a
  finite dimensional linear space.
\end{pr}

\begin{proof}
  By Proposition~\ref{u4}, a complete list of transitive systems of
  three subspaces has the following form:
  $$
  \begin{array}{c}
    \begin{array}{cc}
      S_1=(\mathbb C;0,0,0),&  S_5=(\mathbb C;0,\mathbb C,\mathbb C),\\[1.5mm]
      S_2=(\mathbb C;\mathbb C,0,0),& S_6=(\mathbb C;\mathbb C,0,\mathbb
      C),\\[1.5mm]
      S_3=(\mathbb C;0,\mathbb C,0),& S_7=(\mathbb C;\mathbb C,\mathbb
      C,0),\\[1.5mm]
      S_4=(\mathbb C;0,0,\mathbb C),& S_8=(\mathbb C;\mathbb C,\mathbb
      C,\mathbb C),\\[1.5mm]
    \end{array}
    \\
    S_9=(\mathbb C^2;\mathbb C(1,0),\mathbb C(0,1),\mathbb C(1,1)).
  \end{array}
  $$
  By the result of Section~3, $\Sigma_3=\{0,1,
  {\textstyle\frac32},2,3\}$.
  
  If $\alpha=0$, the only representation of the algebra $\EuScript
  P_{3,0}$, up to equivalence, is $P_1=0$, $P_2=0$, $P_3=0$ on $U=\mathbb
  C$ and, consequently, the system of there subspaces generated by
  this representation is isomorphic to $S_1$.
  
  If $\alpha=1$ there are three inequivalent irreducible
  representations of the algebra $\EuScript P_{3,1}$. The first
  representation is $P_1=I$, $P_2=0$, $P_3=0$ on the space $U=\mathbb
  C$. The system of three subspaces corresponding to this
  representation is isomorphic to $S_2$. The second representation is
  given by $P_1=0$, $P_2=I$, $P_3=0$ on the space $U=\mathbb C$. The
  corresponding system of three subspaces is isomorphic to $S_3$. The
  third representation is $P_1=0$, $P_2=0$, $P_3=I$ on $U=\mathbb
  C$. The corresponding system of three subspaces is isomorphic to $S_4$.

  If $\alpha=3/2$, there is a unique irreducible representation of the
  algebra $\EuScript P_{3,3/2}$. It is given by
  $$
  P_1=\begin{pmatrix}1&0\\0&0\end{pmatrix},\quad
  P_2=\begin{pmatrix}1/4&\sqrt{3}/4\\\sqrt{3}/4&3/4\end{pmatrix},\quad
  P_3=\begin{pmatrix}1/4&-\sqrt{3}/4\\-\sqrt{3}/4&3/4\end{pmatrix}
  $$
  on $U=\mathbb C^2$. The system of three subspaces, corresponding to
  this representation, is transitive and is isomorphic to $S_9$, as
  follows from the complete list in Proposition~\ref{u4} for a finite
  dimensional space.

  If $\alpha=2$, there are three inequivalent irreducible
  representations of $\EuScript P_{3,2}$. For the first
  representation, $P_1=0$, $P_2=I$, $P_3=I$ on the space $U=\mathbb
  C$, the system of subspaces is isomorphic to $S_5$. For the second
  representation, $P_1=I$, $P_2=0$, $P_3=I$ on $U=\mathbb C$, the
  corresponding system is isomorphic to $S_6$. The third representation
  is given by $P_1=I$, $P_2=I$, $P_3=0$ on the space $U=\mathbb
  C$. The system of three subspaces, generated by this representation,
  is isomorphic to $S_7$.
  
  For $\alpha=3$, the unique irreducible representation of $\EuScript
  P_{3,3}$, up to equivalence, is $P_1=I$, $P_2=I$, $P_3=I$ on the
  space $U=\mathbb
  C$ and, hence, the corresponding system of three subspaces is
  isomorphic to $S_8$.
\end{proof}

\subsection{Transitive systems of subspaces, generated by $\Rep\EuScript
P_{4,com}$}

An important tool used for describing the set $\Sigma_n$ for $n\geq4$
and constructing the representations, $\Rep\EuScript P_{4,\alpha}$,
that generate systems of the subspaces $S(u,\pm1)$, $u=3,4,5,\ldots$,
and $S(u,\pm2)$, $u=3,5,7,\ldots$, in the list (\ref{e3}) are the
Coxeter functors, which were constructed in~\cite{KRS}, between the
categories of $*$-representations of $\EuScript P_{n,\alpha}$ for
different values of the parameters.

Let us define a functor $\EuScript T:\Rep\EuScript
P_{n,\alpha}\rightarrow \Rep\EuScript P_{n,n-\alpha}$, which is the
first functor constructed in~\cite{KRS}. Let the orthogonal
projections $P_1$, $P_2$, $\ldots$, $P_n$ be a representation in
$\Rep\EuScript P_{n,\alpha}$ with the representation space $H$. Then
the orthogonal projections $I-P_1$, $I-P_2$, $\ldots$, $I-P_n$
constitute a representation in $\EuScript T(\Rep\EuScript
P_{n,\alpha})$ with the same representation space. The second functor
in~\cite{KRS}, $\EuScript S:\Rep\EuScript P_{n,\alpha}\rightarrow
\Rep\EuScript P_{n,\frac{\alpha}{\alpha-1}}$, is defined as follows.
Again denote by $P_1$, $P_2$, $\ldots$, $P_n$ the orthogonal
projections in $\Rep\EuScript P_{n,\alpha}$ with the representation
space $H$. Let $\Gamma_k:\Imp P_k\rightarrow H$, $k=\overline{1,n}$,
be the natural isometries and
$\Gamma=[\Gamma_1,\Gamma_2,\ldots,\Gamma_n]:\EuScript H=\Imp
P_1\oplus\Imp P_2\oplus\ldots\Imp P_n\rightarrow H$. Then the natural
isometry $\sqrt{\frac{\alpha-1}{\alpha}}\Delta^*$ from the orthogonal
complement in $\hat{H}$ to the subspace $\Imp \Gamma^*$ in $\EuScript
H$ gives the isometries $\Delta_k=\Delta|_{\Imp P_k}:\Imp
P_k\rightarrow\hat{H}$, $k=\overline{1,n}$. The orthogonal projections
$Q_k=\Delta_k\Delta_k^*$, $k=\overline{1,n}$, on the space $\hat{H}$
give the corresponding representation in $\EuScript S(\Rep\EuScript
P_{n,\alpha})$.

\begin{lem}\label{l1}
  The functors $\EuScript T$ and $\EuScript S$ take representations
  that define transitive systems into representations that generate
  transitive systems.
\end{lem}

\begin{proof}
  Property~\ref{p2} immediately proves the statement for the functor
  $\EuScript T$.
  
  Consider now the functor $\EuScript S$. Let a collection of
  orthogonal projections $P_1$, $P_2$, $\ldots$, $P_n$ on a Hilbert
  space $H$ satisfy the condition $\sum_{i=1}^nP_i =\alpha I_H$ for
  some $\alpha$, and the corresponding system of subspaces be
  transitive. Consider the representation $Q_1$, $Q_2$, $\ldots$,
  $Q_n$, $\sum_{k=1}^nQ_k=\frac{\alpha}{\alpha-1}I_{\hat{H}}$, with
  the representation space $\hat{H}$, into which the functor
  $\EuScript S$ maps the representation $P_1$, $P_2$,
  $\ldots$, $P_n$. Let us prove that the system of subspaces generated
  by the representation $Q_1$, $Q_2$, $\ldots$, $Q_n$, that is, the
  system $\hat{S}=(\hat{H};Q_1\hat{H},Q_2\hat{H}, \ldots,Q_n\hat{H})$
  is transitive. Let $R\in \End(\hat{S})$. Then
  \begin{equation}\label{e4}
    Q_kRQ_k=RQ_k,\quad \forall k=\overline{1,n}.
  \end{equation}
  
  Denote by $\hat{C}$ the operator such that
  $\hat{C}:\hat{H}\rightarrow\hat{H}$ and $\hat{C}^*=R$. It follows
  from~(\ref{e4}) that
  \begin{equation}\label{e5}
    Q_k\hat{C}Q_k=Q_k\hat{C},\quad \forall k=\overline{1,n}.
  \end{equation}
  Consider the operators $C_k:\Imp P_k\rightarrow\Imp P_k$,
  $(k=\overline{1,n})$, given by
  \begin{equation}\label{e6}
    C_k=\Delta^*_k\hat{C}\Delta_k,\quad k=\overline{1,n},
  \end{equation}
  and show that the operator $\hat{C}$ can be represented as
  \begin{equation}\label{e7}
    \hat{C}=\frac{\alpha-1}{\alpha}\sum_{k=1}^n\Delta_kC_k\Delta^*_k.
  \end{equation}
  Indeed, using~(\ref{e6}) and the definition of $Q_k$ we get
  \begin{multline*}
    \frac{\alpha-1}{\alpha}\sum_{k=1}^n\Delta_kC_k\Delta^*_k=
    \frac{\alpha-1}{\alpha}\sum_{k=1}^n\Delta_k\Delta^*_k\hat{C}\Delta_k
    \Delta^*_k=\frac{\alpha-1}{\alpha}\sum_{k=1}^nQ_k\hat{C}Q_k=\\
    =\frac{\alpha-1}{\alpha}\sum_{k=1}^nQ_k\hat{C}=
    \frac{\alpha-1}{\alpha}(\sum_{k=1}^nQ_k)\hat{C}=\hat{C}.
  \end{multline*}
  Now, (\ref{e5}) and (\ref{e6}) yield 
  \begin{equation}\label{e8}
    \Delta^*_k\hat{C}=C_k\Delta^*_k,\quad \forall k=\overline{1,n},
  \end{equation}
  and
  \begin{multline*}
    C_k\Delta^*_k=(\Delta^*_k\hat{C}\Delta_k)\Delta^*_k=
    \Delta^*_k\hat{C}(\Delta_k\Delta^*_k)
    =\Delta^*_k\hat{C}Q_k=I_{\Imp P_k}\Delta^*_k\hat{C}Q_k=\\
    (\Delta^*_k\Delta_k)\Delta^*_k\hat{C}Q_k=\Delta^*_k(\Delta_k\Delta^*_k)\hat{C}Q_k
    =\Delta^*_kQ_k\hat{C}Q_k=\Delta^*_kQ_k\hat{C}=\qquad\\
    =\Delta^*_k(\Delta_k\Delta^*_k)\hat{C}=
    (\Delta^*_k\Delta_k)\Delta^*_k\hat{C}=\Delta^*_k\hat{C}.
  \end{multline*}  
  
Consider the operator 
  \begin{equation}\label{e9}
    C=\frac1{\alpha}\sum_{i=1}^n\Gamma_iC_i\Gamma^*_i.
  \end{equation}
  Using properties of the operators $\{\Gamma_i\}_{i=1}^n$,
  $\{\Gamma_i^*\}_{i=1}^n$, $\{\Delta_i\}_{i=1}^n$,
  $\{\Delta_i^*\}_{i=1}^n$,
  \begin{equation}\label{e10}
    \sum_{i=1}^n\Gamma_i\Delta_i^*=0,
  \end{equation}
  \begin{equation}\label{e11}
    \Gamma_i^*\Gamma_j=-(\alpha-1)\Delta_i^*\Delta_j,\quad\mbox{для}\;
    i\neq j,
  \end{equation}
  it follows from~\cite{KRS} that
  \begin{equation}\label{e12}
    C\Gamma_k=\Gamma_kC_k\quad\forall k=\overline{1,n},
  \end{equation}
  \begin{equation}\label{e13}
    C_k=\Gamma_k^*C\Gamma_k\quad\forall k=\overline{1,n},
  \end{equation}
  Indeed,
  \begin{multline*}
    C\Gamma_k=\frac1{\alpha}\sum_{i=1}^n\Gamma_iC_i\Gamma_i^*\Gamma_k=
    \frac1{\alpha}\Gamma_kC_k+\frac1{\alpha}\sum_{\substack{i=1\\ i\neq j}
    }^n\Gamma_iC_i(\Gamma_i^*\Gamma_k)
    =\frac1{\alpha}\Gamma_kC_k-\\\frac{\alpha-1}{\alpha}\sum_{\substack{i=1\\
      i\neq j} }^n\Gamma_i(C_i\Delta_i^*)\Delta_k=
    \frac1{\alpha}\Gamma_kC_k-\frac{\alpha-1}{\alpha}\sum_{\substack{i=1\\
      i\neq j}
    }^n\Gamma_i(\Delta_i^*\hat{C})\Delta_k=\frac1{\alpha}\Gamma_kC_k+\\
    +\frac{\alpha-1}{\alpha}\Gamma_k\Delta_k^*\hat{C}\Delta_k=\Gamma_kC_k
  \end{multline*}
  and
  \begin{multline*}
    \Gamma_k^*C\Gamma_k=\frac1{\alpha}\Gamma_k^*(\sum_{i=1}^n\Gamma_iC_i\Gamma_i^*)\Gamma_k=
    \frac1{\alpha}C_k+\frac1{\alpha}\sum_{\substack{i=1\\ i\neq
      j}}^n\Gamma_k^*\Gamma_iC_i\Gamma_i^*\Gamma_k=
    \frac1{\alpha}C_k+\\+\frac{(\alpha-1)^2}{\alpha}\sum_{\substack{i=1\\
      i\neq
      j}}^n\Delta_k^*\Delta_iC_i\Delta_i^*\Delta_k=\frac1{\alpha}C_k+(\alpha-1)\Delta_k^*\hat{C}\Delta_k
    -\frac{(\alpha-1)^2}{\alpha}C_k=C_k.
  \end{multline*}
  
  It follows from~(\ref{e12}), (\ref{e13}) that
  $CP_k=C\Gamma_k\Gamma_k^*=\Gamma_kC_k\Gamma_k^*=
  \Gamma_k\Gamma_k^*C_k\Gamma_k\Gamma_k^*= P_kCP_k$, which means that
  $C\in\End(S)$, where $S=(H;P_1H,P_2H,\ldots,P_nH)$. Because, by the
  assumption, the system $S$ is transitive, we have $\End(S)=\mathbb C
  I_H$ and, consequently, $C$ is a scalar operator. By~(\ref{e13}),
  $C_k=\lambda I_{\Imp P_k}$ $(k=\overline{1,n})$. Now, according
  to~(\ref{e7}),  $\hat{C}=\lambda I_{\hat{H}}$ and, correspondingly,
  $R$ is a scalar operator. This ends the proof.
\end{proof}

\begin{lem}\label{l2}
  The mapping
  \begin{equation}\label{e}
    \lambda=\frac{b^2-a^2c^2}{(1-a^2)^2}+i\frac{2abc}{(1-a^2)^2}
  \end{equation}
  realizes a one-to-one correspondence between the region $\Omega$ and
  the complex plain with the deleted points $0$ and $1$.
\end{lem}
\begin{proof}
  Consider the points $A(1,0,0)$, $B(0,1,0)$, and $C(0,0,1)$ as in
  Fig.~\ref{fig1}. The point $C$ of the unit sphere, which does not
  belong to the region $\Omega$, is mapped by~(\ref{e}) into the
  deleted point $0$ of the complex plain $(\lambda)$, see
  Fig.~\ref{fig2}. The point $B$ of the unite sphere does not belong
  to the region $\Omega$ and is mapped by~(\ref{e4}) into the removed
  point $1$. The points of the arc $CB$, which belong to the region
  $\Omega$, that is, all the points of the arc except for the points
  $C$ and $B$, are mapped by~(\ref{e4}) in a one-to-one manner, into
  points of the interval $(0,1)$ of the real axis.

  \begin{figure}[htb]
    \includegraphics[height=5cm]{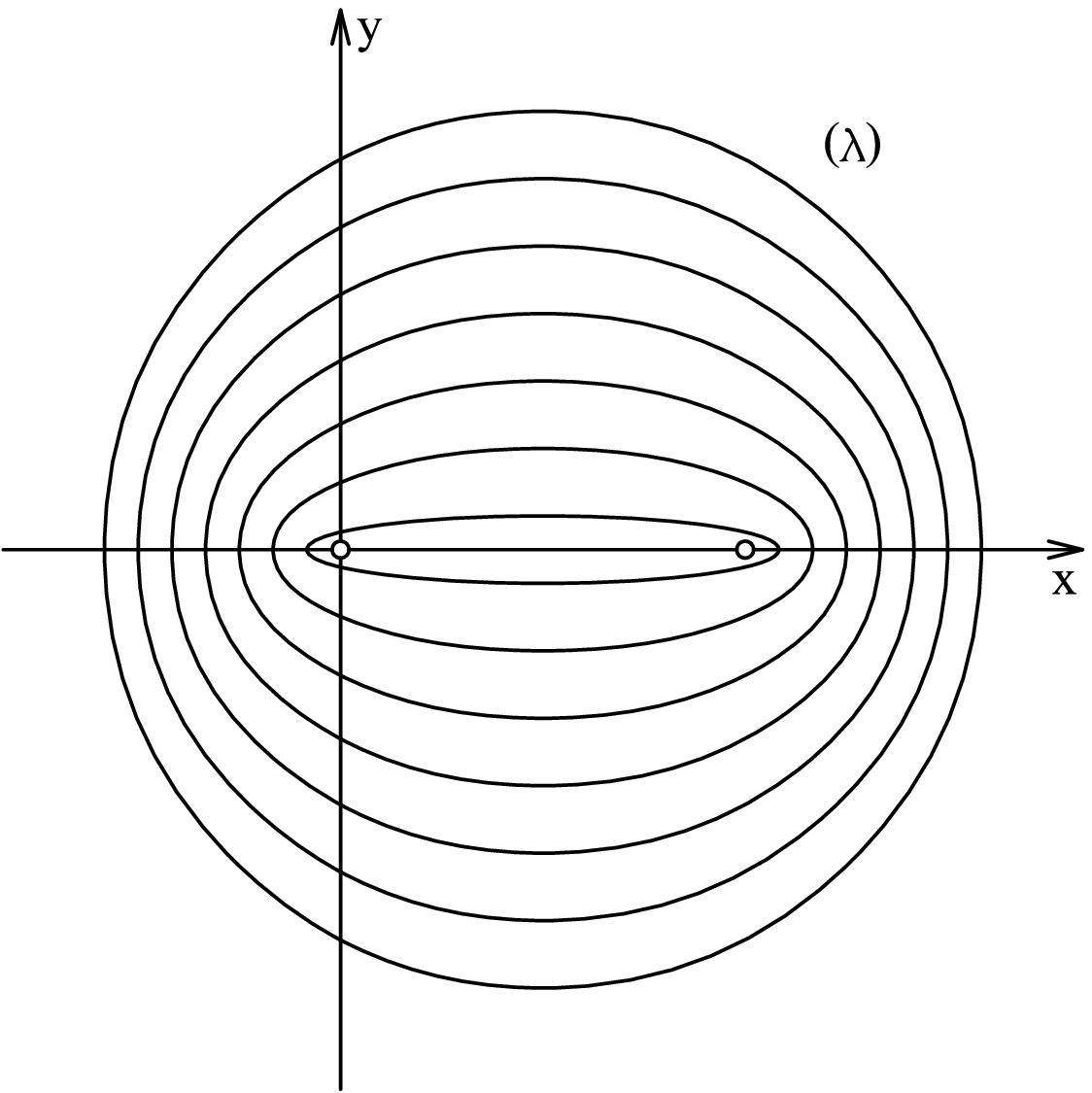}
    \caption{}\label{fig2}
  \end{figure}
  
  Let us fix $0<a<1$. Then $\Gamma_a=\{(a,b,c)\in\mathbb R^3|
  b=\sqrt{1-a^2}\cos \varphi, c=\sqrt{1-a^2}\sin \varphi,
  \varphi\in(-\pi/2,\pi/2] \}\subset\Omega $. Denote
  $\tilde{a}=\frac12\frac{1+a^2}{1-a^2}$ and
  $\tilde{b}=\frac{a}{1-a^2}$. For $x=\Re\lambda$ and $y=\Im\lambda$,
  we get
  $$
  \frac{(x-1/2)^2}{\tilde{a}^2}+\frac{y^2}{\tilde{b}^2}=1,
  $$
  so that the mapping~(\ref{e4}) takes points of the arc $\Gamma_a$,
  in a one-to-one manner, into an ellipse with center in the point
  $(1/2,0)$, major semiaxis $\tilde{a}$ and minor semiaxis $\tilde{b}$.

  As $a\in(0,1)$ ranges from zero to one, the major semiaxis is a
  strictly increasing function with values in the interval
  $(1/2,\infty)$. The minor semiaxis is also a strictly increasing
  function on the interval $(0,1)$ with values $\tilde{b}$ ranging
  over the interval $(0,\infty)$.
\end{proof}

\begin{thm}\label{t8}
  Irreducible nonequivalent $*$-representations of $\EuScript
  P_{4,com}$ generate all transitive systems of four subspaces of a
  finite dimensional linear space.
\end{thm}

\begin{proof}
  By Theorem~\ref{t1}, a complete list of nonisomorphic transitive
  systems of four distinct proper subspaces of a finite dimensional
  linear space is the following:
  $$
  \begin{array}{c}
    B(2,0;\lambda), \quad\lambda\in \mathbb C, \lambda\neq0,1,\\
    B(u,\pm1), \quad u=3,4,5,\ldots,\\
    B(u,\pm2), \quad u=3,5,7,\ldots.
  \end{array}
  $$
  Let us show that the systems $S(2,0;a,b,c)$ are isomorphic to the
  systems $B(2,0;\lambda)$ for
  $\lambda=\frac{b^2-a^2c^2}{(1-a^2)^2}+i\frac{2abc}{(1-a^2)^2}$, up
  to a rearrangement of the subspaces.  Denote $A=1+a$ and $B=b-ic$.
  Then
  $$
  S(2,0;a,b,c)=(\mathbb C^2;\Imp P_1,\Imp P_2,\Imp P_3,\Imp P_4),
  $$
  where
  $$
  \begin{array}{cc}
    \Imp P_1=\mathbb C(A,-B),& \Imp P_3=\mathbb C(B,A),\\
    \Imp P_2=\mathbb C(B,-A),& \Imp P_4=\mathbb C(A,B).
  \end{array}
  $$
  Denote by $R\in M_2(\mathbb C)$ a linear transformation from
  $\mathbb C^2$ to $\mathbb C^2$, such that $R(\Imp P_1)\subset K_1$,
  $R(\Imp P_2)\subset K_2$, $R(\Imp P_4)\subset K_3$, $R(\Imp
  P_3)\subset K_4$. The first three conditions give
  $$
  R=\begin{pmatrix}1& \frac{B}{A}\\
    \frac{A^2+B^2}{2A^2}& \frac{A^2+B^2}{2AB}\end{pmatrix}.
  $$
  The matrix $R$ satisfies the condition $R(\Imp P_3)\subset K_4$
  for $\lambda=\frac{b^2-a^2c^2}{(1-a^2)^2}+i\frac{2abc}{(1-a^2)^2}$.
  In virtue of Lemma~\ref{l2}, this gives an isomorphism, up to a
  rearrangement of the subspaces, between the systems $S(2,0;a,b,c)$,
  where $(a,b,c)\in\Omega$, and the systems $B(2,0;\lambda)$, where
  $\lambda\in \mathbb C$, $\lambda\neq0,1$, for
  $\lambda=\frac{b^2-a^2c^2}{(1-a^2)^2}+i\frac{2abc}{(1-a^2)^2}$. This
  shows that systems that correspond to nonequivalent irreducible
  two-dimensional representations in $\Rep\EuScript P_{4,2}$ are
  nonisomorphic and transitive.

  By Lemma~\ref{l1}, we obtain transitivity, since the dimensions of
  the nonisomorphic systems
  $$
  \begin{array}{l}
    S(u,\pm1), \quad u=3,4,5,\ldots,\\
    S(u,\pm2), \quad u=3,5,7,\ldots,
  \end{array}
  $$
  are different. Since the list of transitive systems, given in
  Section~2, is complete, we have
  $$
  \begin{array}{c}
    S(u,\pm1)\cong B(u,\pm1), \quad u=3,4,5,\ldots,\\
    S(u,\pm2)\cong B(u,\pm2), \quad u=3,5,7,\ldots,
  \end{array}
  $$
  up to a rearrangement of the subspaces.
\end{proof}

In confirmation of the hypothesis formulated in Introduction,
Lemma~\ref{l1} allows to conclude that the system of subspaces,
generated by irreducible $*$-representations of $\EuScript P_{n,com}$
for $n\geq5$ and
$\alpha\in\{\Lambda_n^0,\Lambda_n^1,n-\Lambda_n^1,n-\Lambda_n^0\}$, is
transitive.

\end{document}